\documentclass[conference,letterpaper]{IEEEtran}

\usepackage[margin=1in]{geometry}

\usepackage{amssymb, mathrsfs,verbatim} 
\usepackage[cmex10]{amsmath}
\RequirePackage[numbers]{natbib}
\RequirePackage[colorlinks,linkcolor=black,citecolor=black,urlcolor=blue]{hyperref}
\usepackage{enumerate}
\usepackage[dvips]{graphicx}
\usepackage{verbatim}
\usepackage{bbm}
\usepackage{epsfig}
\usepackage{epstopdf}
\usepackage{graphicx} 

\usepackage{caption}
\usepackage{float,latexsym,amsfonts,amstext,times}
\usepackage{chngcntr}
\usepackage{apptools}
\usepackage{subfig}

\usepackage{color}

\newtheorem{theorem}{Theorem}[section]

\newcommand{\cS}{\mathcal{S}}

\newcommand{\cC}{\mathcal{C}}
\newcommand{\cF}{\mathcal{F}}

\newcommand{\Pro}{\mathsf{P}}
\newcommand{\Exp}{\mathsf{E}}

\newcommand{\bR}{\mathbb{R}}

\DeclareMathOperator\erf{erf}
\begin{document}

\title{Multiple sequences Prophet Inequality Under Observation Constraints.}

\author{\IEEEauthorblockN{
		Aristomenis Tsopelakos, Olgica Milenkovic.}
	\IEEEauthorblockA{Coordinated Science Laboratory, University of Illinois at Urbana-Champaign.}	
}
\vspace{-0.4in}

\maketitle

\begin{abstract}
In our problem, we are given access to a number of sequences of nonnegative i.i.d. random variables, whose realizations are observed sequentially. All sequences are of the same finite length. The goal is to pick one element from each sequence in order to maximize a reward equal to the expected value of the sum of the selections from all sequences. The decision on which element to pick is irrevocable, i.e., rejected observations cannot be revisited. Furthermore, the procedure terminates upon having a single selection from each sequence. Our observation constraint is that we cannot observe the current realization of all sequences at each time instant. Instead, we can observe only a smaller, yet arbitrary, subset of them. Thus, together with a stopping rule that determines whether we choose or reject the sample, the solution requires a sampling rule that determines which sequence to observe at each instant. The problem can be solved via dynamic programming, but with an exponential complexity in the length of the sequences. In order to make the solution computationally tractable, we introduce a decoupling approach and determine each stopping time using either a single-sequence dynamic programming, or a \textit{Prophet Inequality} inspired threshold method, with polynomial complexity in the length of the sequences. We prove that the decoupling approach guarantees at least $0.745$ of the optimal expected reward of the joint problem. In addition, we describe how to efficiently compute the optimal number of samples for each sequence, and its' dependence on the variances.
\end{abstract}

\section{Introduction}
In many applications, multiple data sequences are monitored sequentially with an aim to decide the best instant to terminate the observation procedure, and use the collected information to maximize an objective. Financial applications include the design of posted pricing mechanisms for auctions~\cite{alaei2014bayesian,feldman2014combinatorial,correa2017posted} and contention resolution schemes~\cite{LeeS18}. Recent engineering applications focus on the optimization of computer hardware performance, such as computational sprinting~\cite{huang2019dynasprint,epitropou2019optimal}, which provides a significant performance boost to microchips. These applications, among many others, have motivated a large body of work in optimal stopping theory, with the dynamic programming \cite[Ch. 24]{bhattacharya2021random} being the primary resolution method for a broad class of them.

In systems with multiple data sequences, it is not always possible to observe the current sample of each sequence due to resource limitations or other observation constraints~\cite{nitinawarat2013controlled,xu2021optimum}. For example, when following a large number of auctions in parallel, it is difficult to analyze all offers at all times. In computational sprinting hardware mechanisms, a software predicts the performance boost of a microchip when we allow short-term overheating, but limited computational resources render it impossible to run the software for all microchips, at each instant. In both cases, only a smaller number of sequences can be processed at any given instant, but without restrictions on which they are. This type of observation constraints leads to a problem at the intersection of optimal stopping and multi-armed bandit theory \cite{zhao2022multi}.

In our problem, we assume that the data sequences comprise i.i.d. random variables that are independent of each other, and whose distributions are allowed to differ. Our constraint is that we can only observe a fixed-size subset of the sequences, at each instant. Our goal is to determine which sequences to observe at each instant, and when to stop sampling at each sequence, in order to maximize our reward, which is the sum of the expected values of the observations at the selected stopping times. This is a combined optimal stopping and sampling problem, which could be treated by dynamic programming with an exponential complexity on the length of the sequences. In order to reduce the complexity, we introduce a decoupling approach, based on the Prophet Inequality \cite{correa2017posted,krengel1978semiamarts,assaf2000simple,correa2021prophet,bubna2023prophet}, that reduces the complexity to polynomial, and guarantees at least $0.745$ of the optimal expected reward of the joint problem.

The Prophet Inequality compares the expected value of the picked element at the optimal stopping time, for each sequence, with the expected value of its' maximum element, by providing a tight lower bound on the ratio of the respective expected values. The lower bound is proven to be equal to $0.745$, for an i.i.d. sequence, in \cite{correa2017posted,kertz1986stop}, independent of the distribution and of the length of the sequence. This benchmark has motivated the design of stopping rules for our problem, which always satisfy the Prophet Inequality, although they might be sub-optimal in some cases. In \cite[Corollary 4.7]{correa2017posted}, an algorithm is provided for the computation of the thresholds of such a stopping rule, which is simpler compared to that of the single-sequence dynamic programming (Single-DP), although both methods have linear complexity on the length of the sequences. We refer to the algorithm in \cite{correa2017posted}, as \textit{Prophet Inequality thresholding} method (PI-thresholding).

The basic requirement of the decoupling approach, is the calculation of the optimal number of samples for each sequence, based on which we can apply Single-DP, or PI-thresholding on each sequence separately. This results in a constrained optimization problem that depends on the distributions of all sequences. The computation might be intensive for large sequences, thus, under some smoothness assumptions on the optimization objective, we develop an approximation technique for the number of samples, whose error converges to zero as the length of the sequences goes to infinity.

The paper is organized as follows. In Section II, we present the problem formulation, while in Section III we prove our first result, pertaining to the approximation ratio of the decoupling approach. Then, in Section IV, we describe the approximation technique for the computation of the optimal number of samples for each sequence. In our last Section V, we present computational examples, which support the efficiency of the decoupling approach, for both the Single-DP and the PI-thresholding method. 

\section{Problem formulation} 
Let $(\bR_{+},\cS)$ be an arbitrary measurable space and $(\Omega,\cF,\Pro)$ be a probability space which hosts $M$ independent sequences of $n$ i.i.d. $\bR_{+}$-valued random variables,
\begin{equation*}
X_{i}:= \{X_i(m) \, : \, m \in [n]\}, \quad \quad i \in [M],
\end{equation*}
where $[n]:=\{1,\ldots,n\}$, $[M]:=\{1,\ldots,M\}$. We aim to find the stopping time $\tau_{i} \in [n]$, for each sequence $i \in [M]$, that attains
\begin{equation}\label{opt0}
\sup_{\tau_i \in \mathcal{T}} \Exp[X_{i}(\tau_i)],
\end{equation} 
where $\mathcal{T}$ the class of all \textit{stopping times} that take values in $[n]$. Formally, a stopping time is a random variable $\tau_i$, $i \in [M]$, for which the event $\{ \tau_{i}=m\}$, $m \in [n]$, is fully determined by the observations up to time $m$. Since each sequence $i \in [M]$ is associated to a stopping time $\tau_i$, we make subsequent use of the vector of stopping times of all sequences, i.e., $T := (\tau_1,\ldots,\tau_M)$.

The observations from each sequence $i \in [M]$ are made sequentially and our decision to stop at time $\tau_{i}$ and pick $X_{i}(\tau_i)$ as our ``best" choice for the optimization problem \eqref{opt0} is irrevocable, i.e., we cannot revisit samples we rejected, nor we can examine samples that follow after we stop. The Prophet inequality \cite{correa2017posted,kertz1986stop} provides a tight lower bound of $0.745$ on the ratio of \eqref{opt0} over $\Exp[\max_{m \in [n]} X_{i}(m)]$, for each $i \in [M]$. By the term ``tight", we mean that for any $n$, there exists a distribution that satisfies the inequality by equality, \cite{kertz1986stop}. 

The main aspect that distinguishes our problem from the relevant bibliography, is the introduction of observation constraints to our model, i.e., it is not possible to observe the current element of each sequence at each instant $m \in [n]$, but only from a subset of them of size $K < M$, independent of which they actually are. Hence, we denote by $R(m)$ the subset of sequences that are observed at time $m$, i.e.,
\begin{equation*}
R(m):=\left\{ i \in [M] \, : \, R_{i}(m)=1 \right\},
\end{equation*}
where 
\begin{equation*}
R_{i}(m) := \mathbf{1}\left\{ \mbox{observe} \; X_i(m) \right\}.
\end{equation*}
We say that $R$ is a sampling rule if, for every $m \in [n-1]$, the $R(m+1)$ is $\cF^{R}(m)$-measurable, where $\cF^{R}(m)$ is the $\sigma$-algebra generated by the observed elements up to time $m$ according to rule $R$, i.e.,
\begin{equation*}
\begin{aligned}
\cF^{R}(m):=\qquad \qquad \qquad \qquad \qquad \qquad \qquad \qquad \qquad \\
\begin{cases}
\emptyset, \; &\mbox{ if } \; m=0,\\
\sigma\left(\cF^{R}(m-1), \{X_{i}(m)\, :\, i \in R(m)\}\right), \; &\mbox{ if } \; m \in [n].
\end{cases}
\end{aligned}
\end{equation*}
The policy $(R,T)$ belongs to the class $\cC(K)$ if the number of sequences observed at each sampling instant is equal to $K$, i.e.,
\begin{equation}\label{sampl_cons}
\sum_{i=1}^{M}R_{i}(m) = K, \quad \forall \; m \in [n].
\end{equation}
Our goal is to find a policy $(R,T) \in \cC(K)$ that optimizes the objective  
\begin{equation}\label{st0}
\sup_{(R,T) \in \cC(K)} \sum_{i=1}^{M} \Exp\left[X_i(\tau_i)\right],
\end{equation}
under the assumption that for each $i \in [M]$, the i.i.d. sequence $X_i$ has a finite first moment.

\section{The Decoupling Approach} 
The optimization problem~\eqref{st0} can be solved via dynamic programming, with a computational complexity of $O\left({M \choose K}^n \right)$. In order to reduce the complexity, we describe a decoupling approach that produces a $0.745$-constant approximation for~\eqref{st0}, with complexity polynomial in $n$.

For any sampling rule $R$, we denote by
\begin{equation*}
	N^{R}_{i}(n) := \sum_{m=1}^{n} R_i(m), \quad i \in [M],
\end{equation*}
the total number of elements we have observed from sequence $i$ up to time $n$, and by $\tau^{R}_{i}$ the optimal stopping time of sequence $i$, associated with the sampling rule $R$. The decoupling approach consists of the following steps:
\begin{enumerate}[(i)]
	\item Since for each $i \in [M],$ the sequence $X_i$ is i.i.d., it suffices to determine the number of observations $N_{i}(n)$ for sequence $i$, without pinpointing the exact times at which we observe. The values $\{N_{i}(n) \, :\, i \in [M]\}$ shall satisfy a particular optimization criterion, independent of $R$.	      
	\item We design a sampling rule $R^{d}$ which guarantees $N_{i}(n)$ observations for each $i \in [M]$, and respects the sampling constraint~\eqref{sampl_cons}, which implies 
	\begin{equation}
	\begin{aligned}
	N_{i}^{R^d}(n) &= N_{i}(n), \quad \forall \; i \in [M],\\
	\sum_{i=1}^{M} N_{i}^{R^d}(n) &=K n.
	\end{aligned}
	\end{equation}	      
	\item Given the $N_{i}^{R^d}(n)$ observations for each $i \in [M]$, in order to determine the decoupled optimal stopping times $\tau_{i}^{R^d}$, $i \in [M]$, we can use either the single-sequence dynamic programming~\cite[Chapter 24]{bhattacharya2021random}, or the Prophet Inequality thresholding method in \cite[Corollary 4.7]{correa2017posted}. The latter, although computationally simpler, it may provide a sub-optimal reward. 
\end{enumerate}		
In order to formulate the optimization criterion that the $\{N_{i}(n):i \in [M]\}$ must satisfy, we note that since the elements of each sequence are i.i.d., all subsequences of $N_{i}(n)$ elements have the same statistical behavior, and thus, by Prophet inequality~\cite{correa2017posted, kertz1986stop}, for $R^d$, and for each $i \in [M]$,
\begin{equation}\label{ineq}
\Exp\left[X_i(\tau^{R^d}_i)\right] \geq D\, \Exp\left[ \max_{1 \leq m \leq N^{R^d}_{i}(n)} X_i(m) \right],
\end{equation}
where $D:=0.745$. Thus, the $\{N_{i}(n):i \in [M]\}$ are chosen as the maximizers of
\begin{equation}\label{sxh}
\max_{(n_1,\ldots,n_M)\in \mathfrak{C}(K)} \sum_{i=1}^{M} \Exp\left[ \max_{1 \leq m \leq n_i} X_i(m) \right],
\end{equation}
where
\begin{equation}
\mathfrak{C}(K):=\left\{ (n_1,\ldots,n_M) \in [n]^M \,:\, \sum_{i=1}^{M}=Kn \right\},
\end{equation}
because by \eqref{ineq}, the criterion \eqref{sxh} guarantees that 
\begin{equation}\label{ineqq6}
\begin{aligned}	
&\sup_{(R,T) \in \cC(K)} \sum_{i=1}^{M} \Exp\left[X_i(\tau_i)\right] \\
&\quad\geq D \max_{(n_1,\ldots,n_M)\in \mathfrak{C}(K)} \sum_{i=1}^{M} \Exp\left[ \max_{1 \leq m \leq n_i} X_i(m) \right].
\end{aligned}
\end{equation}

\begin{theorem}
The decoupling approach achieves at least $0.745$ of the optimal expected reward of the joint optimization problem~\eqref{st0}, provided that the number of samples from each sequence is optimized as of \eqref{sxh}.
\end{theorem}

\begin{IEEEproof}
In view of \eqref{ineqq6}, it suffices to show that
\begin{equation}
\begin{aligned}
&\sup_{(R,T) \in \cC(K)} \sum_{i=1}^{M} \Exp\left[X_i(\tau_i)\right] \\
&\quad \leq \max_{(n_1,\ldots,n_M)\in \mathfrak{C}(K)} \sum_{i=1}^{M} \Exp\left[ \max_{1 \leq m \leq n_i} X_i(m) \right].
\end{aligned}
\end{equation}
Let us denote by $R^*$ the optimal sampling rule for \eqref{st0}, as determined by the dynamic programming algorithm, which we denote by $\mathcal{A}$ throughout the proof. Thus,
\begin{equation}
\sup_{(R,T) \in \cC(K)} \sum_{i=1}^{M} \Exp\left[X_i(\tau_i)\right] = \sum_{i=1}^{M} \Exp\left[X_i\left(\tau^{R^*}_i\right)\right].
\end{equation}
The event 
\begin{align}\label{evt}
\{\tau^{R^*}_{1}=t_1,\ldots,\tau^{R^*}_{M}=t_{M}\}
\end{align}
is fully determined by the elements we observed up to the respective times $t_1,\ldots,t_M \in [n]$, i.e.,
\begin{equation}\label{ddt}
\bigcup_{i=1}^{M} \left\{ X_{i}(k)R^*_{i}(k) \, :\, 1 \leq k \leq t_i\right\}.
\end{equation}
The algorithm $\mathcal{A}$, generates also the conditions, which based on the samples of each sequence up to times $\{t_i : i\in [M]\}$ respectively, determine the event \eqref{evt}. Thus, if we denote by $\mathcal{A}_{t_1,..,t_M}$ the conditions of $\mathcal{A}$ for the times $\{t_i : i\in [M]\}$, with a slight abuse of notation we have 
\begin{equation*}
\begin{aligned}
&\Pro\{{\tau^{R^*}_{1}=t_1,..,\tau^{R^*}_{M}=t_{M}\}} \\
&= \Pro\left(\mathcal{A}_{t_1,..,t_M}\left(\bigcup_{i=1}^{M} \left\{ X_{i}(k)R^*_{i}(k) : 1 \leq k \leq t_i\right\}\right)\right).
\end{aligned}
\end{equation*}
Since all sequences are i.i.d. and independent of each other, the random variables in~\eqref{ddt} are interchangeable with the first $N^{R^*}_{i}(t_i)$ random variables in each sequence $i \in [M]$. As a result,
\begin{equation*}
	\begin{aligned}
		&\Pro\left(\tau^{R^*}_{1}=t_1,..,\tau^{R^*}_{M}=t_{M}\right) \\
		&= \Pro\left(\mathcal{A}_{t_1,..,t_M}\left(\bigcup_{i=1}^{M} \left\{ X_{i}(k) : 1 \leq k \leq N^{R^*}_{i}(t_i)\right\}\right)\right).
	\end{aligned}
\end{equation*}
Thus, for a sampling rule $\widetilde{R}$ which examines the elements one by one, without missing the elements that it can not observe, as they remain stand by, one has
\begin{equation}
\sum_{i=1}^{M} \Exp\left[X_i\left(\tau^{R^*}_i\right)\right] = \sum_{i=1}^{M} \Exp\left[ X_i\left(\tau^{\widetilde{R}}_{i}\right) \right],
\end{equation}
and
\begin{equation}
\tau^{\widetilde{R}}_{i} = N^{R^*}_{i}\left(\tau^{R^*}_{i}\right), \quad \mbox{a.s., }\; \forall \; i \in [M].
\end{equation}
Hence, for each $i \in [M]$, 
\begin{equation}
\begin{aligned}	
\Exp\left[ X_i\left(\tau^{\widetilde{R}}_{i}\right) \right] &= \Exp\left[ X_i\left(N^{R^*}_{i}\left(\tau^{R^*}_{i}\right)\right) \right] \\
&\leq \Exp\left[ \max_{1 \leq m \leq N^{R^*}_{i}\left(\tau^{R^*}_{i}\right)} X_i(m) \right].
\end{aligned}
\end{equation}
Therefore, it suffices to show that
\begin{equation}\label{tsh}
\begin{aligned}
&\sum_{i=1}^{M}	\Exp\left[ \max_{1 \leq m \leq N^{R^*}_{i}(\tau^{R^*}_{i})} X_i(m) \right]\\
&\quad \leq \max_{(n_1,\ldots,n_M)\in \mathfrak{C}(K)} \sum_{i=1}^{M} \Exp\left[ \max_{1 \leq m \leq n_i} X_i(m) \right].
\end{aligned}
\end{equation}
Indeed, by sampling constraint~\eqref{sampl_cons}, and since $\tau^{R^*}_{i} \leq n$, for all $i \in [M]$, we have
\begin{equation*}
\sum_{i=1}^{M} N_{i}^{R^*}(\tau^{R^*}_{i}) \leq \sum_{i=1}^{M} \sum_{m=1}^{n} R^{*}_{i}(m)\leq Kn.
\end{equation*}
Since each term in~\eqref{sxh} is increasing in $n_i$, we conclude \eqref{tsh}.
\end{IEEEproof}

\section{Maximization Problem}
We focus on the computation of the maximizers $n_{1}^{*},\ldots,n_{M}^{*}$ of~\eqref{sxh}. For large $n$, finding the exact values is computationally demanding, and thus we suggest an approximation technique along with error guarantees. 

For each $i \in [M]$, we restrict our attention to absolutely continuous random variables, whose probability density and cumulative distribution functions are denoted by $f_i$ and $F_i$, respectively. The density of the maximum of $n_i$ observations from sequence $i$ is denoted by $g_i$, i.e.,
\begin{equation*}
g_{i}(x) := n_i \left( F_{i}(x) \right)^{n_i -1} f_{i}(x), \quad x \geq 0.
\end{equation*}
The maximization problem~\eqref{sxh} turns into \textbf{Problem $\mathbf{P_1}$}:
\begin{equation}\label{of}
\max_{n_1,\ldots,n_M} \sum_{i=1}^{M} n_i \int_{0}^{\infty} x\left( F_i(x) \right)^{n_i-1} f_i(x) dx, 
\end{equation}
subject to
\begin{equation}\label{c1}
\begin{aligned}
\sum_{i=1}^{M} n_i = K n,
\end{aligned}
\end{equation}
and $n \geq n_i \geq 0$, for all  $i \in [M]$.

A solution approach to $\mathbf{P_1}$ is the exhaustive search, of computational complexity $O\left(n^M\right)$. Thus, in the following subsection, we introduce a computationally simpler approximation method, whose computational complexity is $O\left(2^M\right)$, and hence independent of $n$, with an approximation error converging to $0$ as $n \to \infty$.

\subsection{The approximation problem}
For each $i \in [M]$, we set $c_i := n_i/n$ and rewrite problem \textbf{$\mathbf{P_1}$} as \textbf{Problem $\mathbf{P_2}$}: 
\begin{equation}\label{obj_g}
\max_{c_1,\ldots,c_M} \sum_{i=1}^{M} c_i n \int_{0}^{\infty} x\left( F_i(x) \right)^{c_i n-1} f_i(x) dx, 
\end{equation}
subject to
\begin{equation}\label{cst1}
\sum_{i=1}^{M} c_i = K,
\end{equation}
and $1 \geq c_i \geq 0$, for all  $i \in [M]$. We denote by $G(c_1,\ldots,c_M)$ the objective function~\eqref{obj_g}, where $(c_1,\ldots,c_M)$ lies in the set
\begin{equation*}
\mathcal{D} := \left\{ (c_1,\ldots,c_M) \in [0,1]^{M} \, :\, \sum_{i=1}^{M} c_i = K\right\}.
\end{equation*}
For simplicity, we assume that the function $G(c_1,\ldots,c_M)$ is differentiable everywhere on $\mathcal{D}$. Thus, the function $G(c_1,\ldots,c_M)$ is continuous everywhere on the compact set $\mathcal{D}$, and by the Extreme value theorem~\cite[Theorem 4.16]{rudin1953principles} it achieves a maximum on $\mathcal{D}$. By the Lagrange multipliers method and provided that the induced system of equations has a solution, we obtain the solution of \textbf{$\mathbf{P_2}$}, denoted by $(\hat{c}_1,\ldots,\hat{c}_M)$. Then, within  $\mathbf{P_1}$, we replace $n \geq n_i \geq 0$ by
\begin{equation*}
\lfloor \hat{c}_{i} n \rfloor \leq n_i \leq \lceil \hat{c}_i n \rceil, \quad \forall \; i \in [M],
\end{equation*}
which reduces the number of possible solutions from $n+1$ to $2$, for all $i \in [M]$. Hence, we obtain \textbf{Problem $\mathbf{P_3}$}:
\begin{equation}\label{objj}
\max_{n_1,\ldots,n_M} \sum_{i=1}^{M} n_i \int_{0}^{\infty} x\left( F_i(x) \right)^{n_i-1} f_i(x) dx, 
\end{equation}
subject to 
\begin{equation}\label{css2}
\begin{aligned}
\sum_{i=1}^{M} n_i &= K n,\\
\lfloor \hat{c}_i n \rfloor \leq n_i &\leq  \lceil \hat{c}_i n \rceil, \quad \forall \; i \in [M]. 
\end{aligned}
\end{equation}

Solving $\mathbf{P_3}$ using exhaustive search, we compute an approximation of the optimal solution of $\mathbf{P_1}$.

\subsection{The cost of approximation}
We provide a bound on the approximation error for the sample sizes $n_i$, which converges to $0$ as $n$ increases. We denote by $H(n_1,\ldots,n_M)$ the objective function of problems $\mathbf{P_1}$, $\mathbf{P_3}$, and by
\begin{enumerate}[(i)]
	\item $(n^*_1,\ldots,n^*_M)$ the optimal solution of problem $\mathbf{P_1}$, with $c^*_i :=n^*_i/n$, for all $i \in [M]$.
	
	\item $(\check{n}_1,\ldots,\check{n}_M)$ the optimal solution of problem $\mathbf{P_3}$, with $\check{c}_{i} := \check{n}_{i}/n$, for all $ i \in [M]$.
\end{enumerate}

We define the cost of approximation as 
\begin{equation}\label{cost}
C:=\big{|} H(n^*_1,\ldots,n^*_M) - H(\check{n}_1,\ldots,\check{n}_M)\big{|}.
\end{equation}
Next, we prove that under smoothness conditions on $G$, the cost $C$ is bounded by a constant multiplied by $e$,
\begin{equation}\label{err}
e:=\|(\hat{c}_1,\ldots,\hat{c}_M) - (\check{c}_1,\ldots,\check{c}_M)\|.
\end{equation} 

\begin{theorem}\label{prp1}
	If $G$ is everywhere differentiable on $\mathcal{D}$, with bounded derivative, then $C \leq Q \, e$, where
	\begin{equation}\label{q}
	Q:= \max\left\{ \|\nabla G(c_1,..,c_M) \| :(c_1,..,c_M) {\in} \mathcal{D} \right\}.
	\end{equation}
\end{theorem}

\begin{IEEEproof}
By definition of the objective functions $H$, $G$, we observe that for any $n_1,\ldots,n_M \in [n]$ it holds
\begin{equation}
G\left( \frac{n_1}{n},\ldots,\frac{n_M}{n} \right) = H(n_1,\ldots,n_M),
\end{equation}
which implies that
\begin{equation}\label{t_a}
\begin{aligned}
G(c^*_1,..,c^*_M)=&H(n^*_1,..,n^*_M) \\
&\geq H(\check{n}_1,..,\check{n}_M) = G(\check{c}_1,..,\check{c}_M),
\end{aligned}
\end{equation}
where the inequality follows by the optimality of $(n^*_1,\ldots,n^*_M)$ for problem $\mathbf{P_1}$. Also, by the optimality of $(\hat{c}_1,\ldots,\hat{c}_M)$ for problem $\mathbf{P_2}$, we have
\begin{equation}\label{t_b}
G(\hat{c}_1,\ldots,\hat{c}_M) \geq G(c^*_1,\ldots,c^*_M).
\end{equation}
By definition \eqref{cost}, and inequalities \eqref{t_a}-\eqref{t_b} we have 
\begin{equation}
C \leq \big{|}G(\hat{c}_1,\ldots,\hat{c}_M) - G(\check{c}_1,\ldots,\check{c}_M)\big{|}.
\end{equation}
Since $G$ is everywhere differentiable with bounded derivative on $\mathcal{D}$, and $\mathcal{D}$ is convex, by Rademacher's theorem \cite[Theorem 1.41]{weaver2018lipschitz} we deduce that
\begin{equation}\label{t_c}
\begin{aligned}
	\| G&(\hat{c}_1,\ldots,\hat{c}_M)-G(\check{c}_1,\ldots,\check{c}_M) \| \\
	&\leq Q\, \|(\hat{c}_1,\ldots,\hat{c}_M) - (\check{c}_1,\ldots,\check{c}_M)\|,
\end{aligned}
\end{equation}
where $Q$ is as defined in~\eqref{q}. This proves the claim.
\end{IEEEproof}

\textit{Remark}: By the second constraint in~\eqref{css2}, it holds
\begin{equation}
\hat{c}_i  -1/n \leq \check{c}_i \leq \hat{c}_i  +1/n,
\end{equation}
which implies that $e \to 0$ as $n \to \infty$, and by Theorem~\ref{prp1}, we conclude that $C \to 0$ as $n \to \infty$.

\subsection{Uniform densities example}
We consider $M$ sequences of $n$ i.i.d. random variables each, which follow a $U[a_i,b_i]$ uniform distribution, where $0 < a_i < b_i$, for each $i \in [M]$. We can observe only $K$ sequences at each time instant. In this case, problem \textbf{$\mathbf{P_2}$} takes the form
\begin{equation*}
	\max_{c_1,\ldots,c_M} \left\{ \sum_{i=1}^{M} b_i - \frac{b_i - a_i}{c_i n +1} \right\},
\end{equation*} 
subject to its underlying constraints. By the Lagrange multipliers method, for each $i \in [M]$, we have
\begin{equation}\label{hc1}
	\hat{c}_i = \frac{\sqrt{b_i - a_i}}{\left(\sum_{j=1}^{M}\sqrt{b_j - a_j}\right)/M}\left( \frac{K}{M} + \frac{1}{n} \right) - \frac{1}{n}.
\end{equation}
\textit{Remark}: For all $i \in [M]$, the differences $b_i - a_i$ determine the values of $\hat{c}_i$, which implies that the variances $(b_i -a_i)^2/12,$ govern the sampling rates.

\subsection{Gaussian densities example}
In various frameworks, e.g. \cite[Section 5.1]{huang2019dynasprint}, the sequences follow a Gaussian density $\mathcal{N}(\mu_i,\sigma^2_i)$, for each $i \in [M]$, and we can observe $K$ sequences at each instant. We assume that the means of the Gaussians are large enough, and the variances relatively small, so that the random variables are positive with probability practically equal to one. According to Blom's formula \cite{royston1982algorithm}, and the approximation of the $\erf$ function presented in \cite{dominici2008some}, for $n_i$ large enough, and for each $i \in [M]$,
\begin{equation*}
E\left[\max_{1 \leq m \leq n_i} X_{i}(m) \right]{\simeq} \mu_i +\sigma_i \sqrt{\frac{\pi}{2}} \frac{10}{8} \left(\frac{1}{n_i +0.25} -0.8\right).
\end{equation*}
Thus, problem \textbf{$\mathbf{P_2}$} reduces to
\begin{equation}
\max_{c_1,\ldots,c_M} \sum_{i=1}^{M} \frac{\sigma_i}{c_i \,n +0.25},
\end{equation}
subject to its defining constraints. By Lagrange multipliers method, we obtain
\begin{equation}\label{hc2}
\hat{c}_i = \frac{\sqrt{\sigma_i}}{\left(\sum_{i=1}^{n}\sqrt{\sigma_i} \right)/M} \left(\frac{K}{M}+\frac{0.25}{n}\right), \;\; i \in [M].
\end{equation}
We observe that except for the change in a constant multiplier of $1/n$, the $\hat{c}_i$ in \eqref{hc2} exhibit the same dependence on the variances as those in \eqref{hc1}, especially for large $n$.

\section{Computational examples}
We compare the expected reward of the joint dynamic programming problem \eqref{st0} versus that of the decoupled problem. For the decoupled problem, we precompute the optimal number of observations for each sequence, and then on each individual sequence, we run either (i) a single-sequence dynamic programming (Single-DP), or (ii) a Prophet Inequality thresholding (PI-thresholding)  method \cite[Corollary 4.7]{correa2017posted}.

We consider $M=3$ sequences that follow three different uniform distributions $U[0,3]$, $U[0.5,2.5]$, $U[1,2]$, which have the same mean but different variances, and of lengths $n=5,\ldots, 10$. We consider two cases $K=1$ and $K=2$, and we plot the ratio of the expected reward of the decoupled problem, for both Single-DP and PI-thresholding, over that of the joint problem.  

For $K=1$ the decoupling approach offers a very good approximation for the problem~\eqref{st0}, with a ratio above $0.92$ for the Single-DP case, and above $0.91$ for the PI-thresholding, for all $n$. We also note that the ratio of the PI-thresholding method is no more than $1\%$ smaller than that of the Single-DP method which is always optimal for a single sequence. 

For $K=2$ the ratio is smaller compared to the former case, but it is still a good approximation, with a ratio above $0.88$ for the Single-DP, and above $0.87$ for the PI-thresholding method, for all $n$. In this case the gap between Single-DP and PI-thresholding is larger, but no more than $10\%$.
\begin{figure}[h!]
		\subfloat[$M=3$, $K=1$]{
			\includegraphics[width=0.5\linewidth]{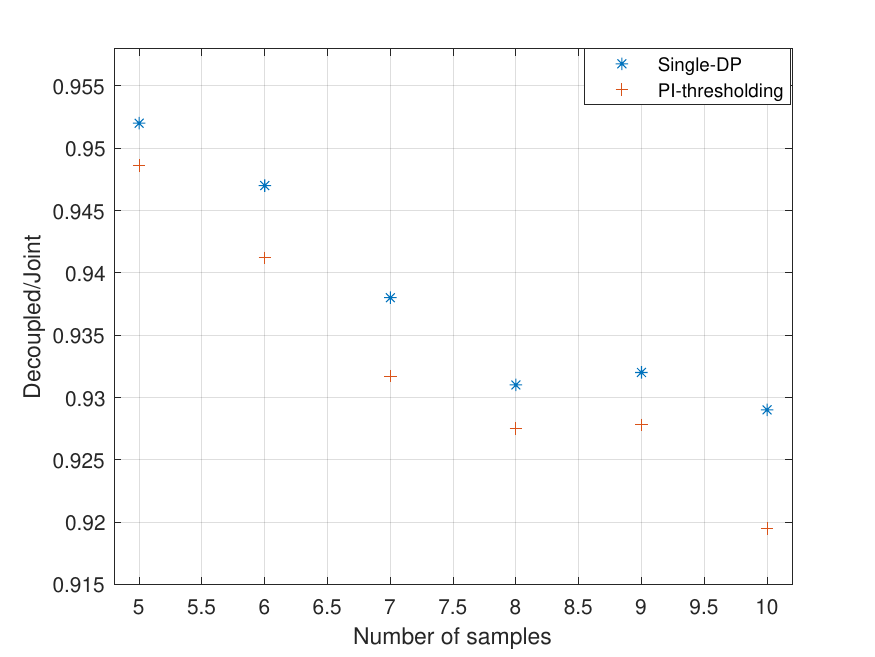}  
			\label{fig:luii}
		}
		\subfloat[$M=3$, $K=2$]{
			\includegraphics[width=0.5\linewidth]{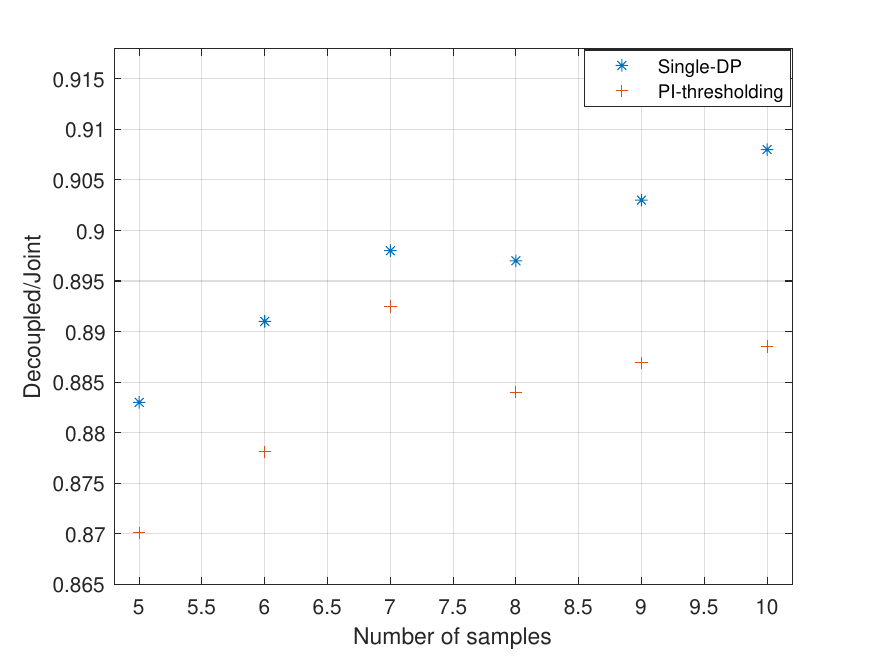}  
			\label{fig:m_ni}
		} 
		\caption{Joint/Decoupled for Single-DP, PI-thresholding.}
		\label{fig:nonhomog_asymi} 
		\vspace{-0.17in}
	\end{figure}
	
\section*{Acknowledgments}
The work was supported in part by NSF awards 2008125 and 1956384, through Coordinated Science Laboratory, University of Illinois at Urbana-Champaign. 

\bibliographystyle{IEEEtran}
\bibliography{biblio_n}

\end{document}